\documentclass[11pt,draft]{amsart}
\usepackage[cp850,latin1]{inputenc}
\usepackage{graphics}

\newtheorem{theorem}{Theorem}
\newtheorem{lemma}[theorem]{Lemma}
\newtheorem{proposition}[theorem]{Proposition}
\newtheorem{corollary}[theorem]{Corollary}

\theoremstyle{definition}

\newtheorem{example}[theorem]{Example}

\theoremstyle{remark}

\usepackage{amssymb,amsmath}

\def\div{\mathop\mathrm{div}\nolimits}

\def\tr{\mathop\mathrm{trace}\nolimits}
\def\det{\mathop\mathrm{det}\nolimits}
\def\Ric{\mathop\mathrm{Ric}\nolimits}
\def\H{\mathbb{H}}

\def\n{\nabla}

\newcommand{\mdos}{\mbox{$\mathbb{M}^2$}}
\newcommand{\mene}{\mbox{$\mathbb{M}^n$}}

\newcommand{\s}{\mbox{$\Sigma$}}

\newcommand{\R}{\mbox{${\mathbb R}$}}

\newcommand{\N}{\mbox{$\mathbb{N}^{n+1}$}}
\newcommand{\g}[2]{\mbox{$\langle #1 ,#2 \rangle$}}
\newcommand{\fle}{\mbox{$\rightarrow$}}
\newcommand{\rf}[1]{\mbox{(\ref{#1})}}
\newcommand{\rl}[1]{{~\ref{#1}}}
\newcommand{\nablabar}{\mbox{$\overline{\nabla}$}}
\newcommand{\fs}{\mbox{$\mathcal{C}^\infty(\s)$}}

\def\beq{\begin{equation}}
\def\eeq{\end{equation}}

\begin{document}

\title[On the scalar curvature of hypersurfaces in spaces with a Killing field]
{On the scalar curvature of hypersurfaces in spaces with a Killing field}

\author{Alma L. Albujer}
\address{Departamento de Matem\'{a}ticas, Universidad de Murcia, E-30100 Espinardo, Murcia, Spain}
\email{albujer@um.es}
\thanks{A.L. Albujer was supported by FPU Grant AP2004-4087 from Secretar\'{\i}a de Estado de Universidades e
Investigaci\'{o}n, MEC Spain.}

\author{Juan A. Aledo}
\address{Departamento de Matem\'{a}ticas, Universidad de Castilla La Mancha,  E-02071 Albacete, Spain}
\email{juanangel.aledo@uclm.es}
\thanks{J.A. Aledo was partially supported by MEC project MTM2007-65249 and Junta de Comunidades de Castilla-La Mancha project PCI-08-0023, Spain.}

\author{Luis J. Al\'\i as}
\address{Departamento de Matem\'{a}ticas, Universidad de Murcia, E-30100 Espinardo, Murcia, Spain}
\email{ljalias@um.es}
\thanks{A.L. Albujer and L.J. Al\'\i as were partially supported by MEC project MTM2007-64504, and Fundaci\'{o}n S\'{e}neca project 04540/GERM/06, Spain.}
\thanks{This research is a result of the activity developed within the framework of the
Programme in Support of Excellence Groups of the Regi\'{o}n de Murcia, Spain, by
Fundaci\'{o}n S\'{e}neca, Regional Agency for Science and Technology (Regional Plan for
Science and Technology 2007-2010).}

\date{April 2008}



\begin{abstract}
We consider compact hypersurfaces in an $(n+1)$-di\-men\-sio\-nal either Riemannian or Lorentzian space \N\ endowed with a conformal Killing vector field. For such hypersurfaces, we establish an integral formula which, especially in
the simpler case when $\N=\mene \times\R$ is a product space, allows us to derive some interesting consequences in terms of the scalar curvature of the hypersurface. For instance, when $n=2$ and \mdos\ is either the sphere
$\mathbb{S}^2$ or the real projective plane $\mathbb{RP}^2$, we characterize the slices of the trivial totally geodesic foliation $\mdos\times\{t\}$ as the only compact two-sided surfaces with constant Gaussian curvature in the
Riemannian product $\mdos\times\mathbb{R}$ such that its angle function does not change sign. When $n\geq 3$ and \mene\ is a compact Einstein Riemannian manifold with positive scalar curvature, we also characterize the slices as the
only compact two-sided hypersurfaces with constant scalar curvature  in the Riemannian product $\mene\times\mathbb{R}$  whose angle function does not change sign. Similar results are also established for spacelike hypersurfaces in a
Lorentzian product $\mathbb{M}\times\mathbb{R}_1$.
\end{abstract}

\maketitle

\section{Introduction}
In the last years, the study of surfaces in product spaces $\mathbb{M}^2\times \mathbb{R}$, where $\mathbb{M}^2$ is a Riemannian surface, has undergone an important development, especially since in the work \cite{AR} Abresch and
Rosenberg generalized Hopf's theorem to homogeneous product spaces. Regarding minimal surfaces in such ambient spaces, they have been of considerable interest since the former works by Rosenberg \cite{Ro} and Meeks and Rosenberg
\cite{MR}. For instance, entire minimal graphs and, more generally, constant mean curvature graphs in $\mathbb{M}^2\times\mathbb{R}$ have been recently studied by several authors, establishing certain Bernstein type results in this
context \cite{Ro,ADR,ER}. Also constant Gaussian curvature surfaces \cite{AEG1} and constant extrinsic curvature surfaces \cite{EGR} in the homogeneous product spaces $\mathbb{H}^2 \times \mathbb{R}$ and $\mathbb{S}^2 \times
\mathbb{R}$ have been recently studied. In fact, in \cite{AEG1} the authors showed that the only complete surfaces of constant Gaussian curvature $K>0$ in $\mathbb{H}^2 \times \mathbb{R}$ (resp. $K>1$ in $\mathbb{S}^2 \times
\mathbb{R}$) are rotational surfaces. In addition, they proved the non-existence of complete surfaces with constant Gaussian curvature $K<-1$ in $\mathbb{H}^2 \times \mathbb{R}$ and $\mathbb{S}^2 \times \mathbb{R}$. On the other hand,
in \cite{EGR} the authors established that the only complete surfaces of constant extrinsic curvature in $\mathbb{S}^2\times\mathbb{R}$ and $\mathbb{H}^2\times\mathbb{R}$ are rotational spheres.

Both $\mathbb{H}^2 \times \mathbb{R}$ and $\mathbb{S}^2 \times\mathbb{R}$, and more generally any product manifold $\mene \times\R$ endowed with the product metric, are manifolds with a globally defined Killing vector field
$\partial_t$. In this work we deal with compact hypersurfaces in an $(n+1)$-dimensional either Riemannian or Lorentzian space \N\, ($n\geq 2$), endowed with a conformal Killing field, which in the case where \N\ is Lorentzian will be
assumed to be timelike. In fact, we obtain an integral formula for such hypersurfaces (Theorem \ref{thintegral}) which, especially in the simpler case when $\N=\mene \times\R$ is a product space, allows us to develop a nice study in
terms of $S_\mathbb{M}$ and $S$, the scalar curvatures of $\mathbb{M}$ and the hypersurface, respectively. More specifically, we are able to characterize the slices of such product spaces from assumptions on $S_\mathbb{M}$ and $S$
(Subsections \ref{SS41} and \ref{SS42}), which become quite powerful either when $n=2$ or when $n\geq 3$ and $\mene$ is Einstein with non-zero scalar curvature.

For instance, when $n=2$ and \mdos\ is either the sphere $\mathbb{S}^2$ or the real projective plane $\mathbb{RP}^2$, we prove that the slices are the only compact two-sided surfaces in the Riemannian product $\mdos\times\mathbb{R}$
whose angle function $\Theta$ does not change sign and have constant Gaussian curvature (Corollary\rl{corriem}). It is worth pointing out that the condition that the angle function does not change sign is the natural generalization to
our ambient spaces of the condition that the image of the Gauss map of a surface in $\mathbb{R}^3$ is contained in a hemisphere, and it is a natural condition to consider if one wants to conclude that the surface is a slice
\cite{AD,ADR,ER}. When \mdos\ is either the sphere $\mathbb{S}^2$ or the real projective plane $\mathbb{RP}^2$, we also prove that the slices are the only compact spacelike surfaces in the Lorentzian product $\mdos\times\mathbb{R}_1$ with constant Gaussian curvature (Corollary\rl{corlor}), without any additional hypothesis. In
general, when $n\geq 3$ we derive analogous results under the assumption that $\mathbb{M}^n$ is Einstein.

Finally, in Section\rl{S5} we apply our results to the case of entire graphs over $\mathbb{M}^2$ with prescribed Gaussian curvature, obtaining several uniqueness results (Corollaries\rl{coro10} and\rl{coro11}). We also give examples
of non-trivial entire graphs over $\mathbb{H}^2$ with constant Gaussian curvature which show the sharpness of our results (Example\rl{exam}).

\section{Preliminaries}
Let \N\ be an $(n+1)$-dimensional either Riemannian or Lorentzian space ($n\geq 2$). In order to simplify our notation, we will denote by $\g{}{}$, without distinction, the corresponding metric tensor, and we will assume that \N\ is
endowed with a conformal Killing field $T \in \mathfrak{X}(\N)$, which in the case where \N\ is Lorentzian will be assumed to be timelike. Recall that a vector field $T$ on \N\ is called a conformal Killing field if the Lie derivative
of the metric tensor $\g{}{}$ with respect to $T$ satisfies $\mathcal{L}_T\g{}{}=2\phi\g{}{}$ for some smooth function $\phi \in \mathcal{C}^\infty(\N)$. Equivalently, \beq \label{conformal}
\g{\nablabar_VT}{W}+\g{V}{\nablabar_WT}=2\phi\g{V}{W} \eeq for all vector fields $V,W \in \mathfrak{X}(\N)$, where $\nablabar$ denotes the Levi-Civita connection in \N. In particular, $T$ is called a homothetic Killing field if the
function $\phi$ is constant, and just a Killing field whenever that constant vanishes.

Let us consider $\Sigma^n$ a connected hypersurface immersed into \N. In the case where \N\ is Lorentzian, we will assume that $\s$ is a spacelike hypersurface, that is, the metric induced on $\s$ via the immersion is a Riemannian
metric. Since $T$ is a globally defined timelike vector field on \N, it follows that there exists a unique unitary timelike normal field $N$ globally defined on \s\ which is in the same time-orientation as $T$, so that
\[
\g{N}{T}\leq -|T|=-\sqrt{-\g{T}{T}}<0 \quad \mathrm{on} \quad \s.
\]
In that case we will refer to $N$ as the future-pointing Gauss map of \s. On the other hand, when \N\ is Riemannian \s\ is assumed to be a \textit{two-sided} hypersurface in \N. This condition means that there is a globally defined
unit normal vector field $N$. In both cases, we will denote by $\Theta:\s \fle \R$ the smooth function on \s\ given by $\Theta=\g{N}{T}$, which in the Lorentzian case will be always negative.

Let $\n$ denotes the Levi-Civita connection on \s. Then the Gauss and Weingarten formulae for the hypersurface in \N\ are given, respectively, by \beq \label{fgauss} \nablabar_XY=\n_XY+\varepsilon\g{AX}{Y}N, \eeq and \beq
\label{fwein} A(X)=-\nablabar_XN, \eeq for all tangent vector fields $X,Y \in \mathfrak{X}(\s)$. Here $A:\mathfrak{X}(\s) \fle \mathfrak{X}(\s)$ denotes the shape operator (or second fundamental form) of \s\ with respect to $N$, and
$\varepsilon=\g{N}{N}$. As is well known, the curvature tensor $R$ of the hypersurface \s\ is described in terms of the shape operator and the curvature tensor $\overline{R}$ of the ambient space \N\ by the so-called Gauss equation,
which can be written as \beq \label{eqgauss} R(X,Y)Z=(\overline{R}(X,Y)Z)^\top+\varepsilon(\g{AX}{Z}AY - \g{AY}{Z}AX) \eeq for all tangent vector fields $X,Y,Z \in \mathfrak{X}(\s)$, where $(\;)^\top$ denotes the tangential component
of a vector field in $\mathfrak{X}(\N)$ along $\s^n$. Observe that our criterion here for the definition of the curvature tensor is the one in \cite{ONe}. Consider a local orthonormal frame $\{E_1,...E_n\}$ of $\mathfrak{X}(\s)$. We
can also assume that $\{E_1,...,E_n\}$ diagonalizes the shape operator, that is $AE_i=\kappa_iE_i, \, 1\leq i\leq n$ where $\kappa_i,\,1\leq i \leq n$, stand for the principal curvatures of $\s^n$. Then, it follows from the Gauss
equation \rf{eqgauss} that
\begin{eqnarray*}
\Ric(E_j,E_j)&=&\sum_{i \neq j}\g{\overline{R}(E_i,E_j)E_i}{E_j}+\varepsilon \sum_{i \neq j}\kappa_i\kappa_j\\
&=&\overline{\Ric}(E_j,E_j)-\varepsilon\g{\overline{R}(E_j,N)E_j}{N}+\varepsilon \sum_{i \neq j}\kappa_i\kappa_j
\end{eqnarray*}
for $1 \leq j \leq n$, where $\Ric$ and $\overline{\Ric}$ stand for the Ricci operators of $\s^n$ and \N\, respectively. Taking traces, we can also obtain the analogous relation between the scalar curvatures $S$ and $\overline{S}$ of
$\s^n$ and \N, which is given by
\begin{eqnarray}
\label{eqgaussaux}
\nonumber S & = & \sum_{j=1}^n \overline{\Ric}(E_j,E_j)-\varepsilon\overline{\Ric}(N,N)+2\varepsilon\sum_{i<j}\kappa_i\kappa_j\\
{} & = & \overline{S}-2\varepsilon\overline{\Ric}(N,N)+2\varepsilon\sum_{i<j}\kappa_i\kappa_j.
\end{eqnarray}
Recall that the mean curvature function of the hypersurface $\Sigma$ is defined as
\[
H=\frac{\varepsilon}{n}\mathrm{tr}(A)=\frac{\varepsilon}{n}(\kappa_1+...+\kappa_n).
\]
Therefore, the following identity is verified, \beq \label{eq0} n^2H^2-\|A\|^2=2\sum_{i<j}\kappa_i\kappa_j \eeq being $\|A\|^2=\tr{A^2}$ and so \rf{eqgaussaux} becomes
\[
S=\overline{S}-2\varepsilon\overline{\Ric}(N,N)+\varepsilon(n^2H^2-\|A\|^2).
\]
In other words, \beq \label{eqgausssc} \|A\|^2+\overline{\Ric}(N,N)=\varepsilon(\overline{S}-S)-\overline{\Ric}(N,N)+n^2H^2. \eeq

On the other hand, the Codazzi equation of the hypersurface describes the normal component of $\overline{R}(X,Y)Z$ in terms of the derivative of the shape operator, and it is given by \beq \label{eqcodazzi}
\g{\overline{R}(X,Y)Z}{N}=\g{(\n_YA)X-(\n_XA)Y}{Z}, \eeq where $\n_XA$ denotes the covariant derivative of $A$, that is
\[
(\n_XA)Y=\n_X(AY)-A(\n_XY).
\]

\section{An integral formula}
In all what follows, \N\ denotes an $(n+1)$-dimensional either Riemannian or Lorentzian space endowed with a conformal Killing field $T$ (which is assumed to be timelike when \N\ is Lorentzian), and $\Sigma^n$ is an immersed
hypersurface in \N, which is assumed to be spacelike when \N\ is Lorentzian and two-sided when \N\ is Riemannian. Before giving our main result of this section, we need to give an expression for the Laplacian of the function
$\Theta=\g{N}{T}$. In the case where \N\ is a Riemannian manifold, Al\'{\i}as, Dajczer and Ripoll gave in \cite{ADR} an expression for the Laplacian of $\Theta$ in terms of the Ricci tensor of the ambient space and of the norm of the
second fundamental form. Barros, Brasil and Caminha proved in \cite{BBC} the analogous expression when \N\ is Lorentzian. For our purpose, we need this expression in a slightly different way. For that reason, and for the sake of
completeness, we will give here an alternative proof for the expression we will need.

\begin{proposition}\label{prlaplaT}
Under all the assumptions above, the Laplacian of the smooth function $\Theta\in \fs$, $\Theta=\g{N}{T}$, is given by \beq \label{laplatheta} \Delta\Theta=-\varepsilon n \g{\n H}{T}+\Theta (S-\overline{S}+\varepsilon
(\overline{\Ric}(N,N)-n^2H^2))-n(\varepsilon H\phi+\frac{\partial \phi}{\partial N}), \eeq where $\n H$ denotes the gradient of the mean curvature function on \s.
\end{proposition}

\begin{proof}
As $T$ is a conformal Killing vector field, from \rf{conformal} and Weingarten formula \rf{fwein} we have
\[
X(\Theta)=\g{\nablabar_XN}{T}+\g{N}{\nablabar_XT}=-\g{AX}{T}-\g{X}{\nablabar_NT}
\]
for every tangent vector field $X \in \mathfrak{X}(\s)$. Therefore, the gradient of $\Theta$ on \s\ is given by \beq \label{gradtheta} \n \Theta=-A T^\top-(\nablabar_NT)^\top, \eeq and its Laplacian is given by
\beq
\label{laplaT}
\Delta \Theta=-\div(AT^\top)-\div((\nablabar_NT)^\top).
\eeq
We will start by computing $\div(AT^\top)$. By Codazzi equation, it holds that for any tangent vector field $X \in \mathfrak{X}(\s)$
\begin{eqnarray*}
\g{\n_X(AT^\top)}{X}&=&\g{(\n_XA)T^\top}{X}+\g{A(\n_XT^\top)}{X}\\
&=&\g{(\n_{T^\top}A)X}{X}-\g{\overline{R}(T^\top,X)N}{X}+\g{\n_XT^\top}{AX}.
\end{eqnarray*}
Consequently, \beq \label{eq1} \div(AT^\top)=\tr(\n_{T^\top}A)-\overline{\Ric}(T^\top,N)+\sum_{i=1}^n\g{\n_{E_i}T^\top}{AE_i} \eeq where $\{E_1,..,E_n\}$ is again a local orthonormal frame of $\mathfrak{X}(\s)$ which diagonalizes $A$.
From the decomposition $T=T^\top+\varepsilon \Theta N$, we also obtain
\[
\nablabar_XT=\nablabar_XT^\top+\varepsilon X(\Theta)N+\varepsilon \Theta\nablabar_XN,
\]
for every tangent vector field $X \in \mathfrak{X}(\s)$. Then, applying the Weingarten formula \rf{fwein} we get \beq \label{eq2} \n_XT^\top=(\nablabar_XT)^\top+\varepsilon\Theta AX. \eeq Therefore,
\[
\sum_{i=1}^n\g{\n_{E_i}T^\top}{AE_i}=\sum_{i=1}^n\g{\nablabar_{E_i}T}{AE_i}+\varepsilon\Theta \|A\|^2.
\]
Since $\{E_1,...,E_n\}$ diagonalizes $A$, it holds that $\g{\nablabar_{E_i}T}{AE_i}=\g{\nablabar_{AE_i}T}{E_i}$, and so from \rf{conformal}
\[
\g{\nablabar_{E_i}T}{AE_i}=\phi\g{AE_i}{E_i}.
\]

Then, \rf{eq1} becomes
\begin{eqnarray} \label{div1}
\div(AT^\top)&=&\tr(\n_{T^\top}A)+\varepsilon nH\phi+\varepsilon\Theta \|A\|^2-\overline{\Ric}(T^\top,N)\\
\nonumber &=&\varepsilon n \g{\n H}{T}+\varepsilon n H \phi+\varepsilon \Theta \|A\|^2-\overline{\Ric}(T^\top,N),
\end{eqnarray}
where we have used the fact that the trace commutes with the covariant derivative.

It remains to obtain an expression for $\div((\nablabar_NT)^\top)$. Observe that from \rf{conformal} we have
\[
\nablabar_NT=(\nablabar_NT)^\top+\varepsilon\g{\nablabar_NT}{N}N=(\nablabar_NT)^\top+\phi N,
\]
which by the Gauss and Weingarten formulae yields that
\[
\nablabar_X\nablabar_NT=\n_X(\nablabar_NT)^\top+\varepsilon\g{AX}{\n_NT}N+X(\phi)N-\phi AX.
\]
Therefore, \beq \label{eq2.5} \div((\nablabar_NT)^\top)=\sum_{i=1}^n \g{\nablabar_{E_i}(\nablabar_NT)}{E_i}+\varepsilon nH\phi. \eeq Observe that
\[
\overline{R}(E_i,N)T=
-\nablabar_{AE_i}T-\nablabar_{\nablabar_NE_i}T-\nablabar_{E_i}\nablabar_NT+\nablabar_N\nablabar_{E_i}T.
\]
Taking traces in the last expression we easily obtain from \rf{eq2.5}
\begin{eqnarray}\label{eq3}
\hspace*{0.7cm}\div((\nablabar_NT)^\top)&= &\varepsilon n \phi H+\overline{\Ric}(T,N)\\
\nonumber&+&\sum_{i=1}^n\g{\nablabar_N\nablabar_{E_i}T}{E_i}-\sum_{i=1}^n\g{\nablabar_{AE_i}T}{E_i}-\sum_{i=1}^n\g{\nablabar_{\nablabar_NE_i}T}{E_i}.
\end{eqnarray}
To simplify the different terms in \rf{eq3}, applying \rf{conformal} we get the relation
\[
\g{\nablabar_N\nablabar_{E_i}T}{E_i}=N(\g{\nablabar_{E_i}T}{E_i})-\g{\nablabar_{E_i}T}{\nablabar_NE_i}=\frac{\partial \phi}{\partial N}+\g{\nablabar_{\nablabar_NE_i}T}{E_i}.
\]
On the other hand, we have
\[
\sum_{i=1}^n\g{\nablabar_{AE_i}T}{E_i}=\phi \tr(A)=\varepsilon n \phi H.
\]
Turning back to \rf{eq3}, it becomes \beq \label{div2} \div((\nablabar_NT)^\top)=\overline{\Ric}(T,N)+n \frac{\partial \phi}{\partial N}. \eeq Finally, from \rf{eqgausssc}, \rf{laplaT}, \rf{div1} and \rf{div2} and since
$T=T^\top+\varepsilon \Theta N$,
\begin{eqnarray*}
\Delta \Theta&=&-\varepsilon n \g{\n H}{T}-\varepsilon \Theta \|A\|^2+\overline{\Ric}(T^\top,N)-\overline{\Ric}(T,N)-n(\varepsilon H \phi+\frac{\partial \phi}{\partial N})\\
{} & = &-\varepsilon n \g{\n H}{T}-\varepsilon \Theta (\|A\|^2+\overline{\Ric}(N,N))-n(\varepsilon H\phi+\frac{\partial \phi}{\partial N})\\
&=& -\varepsilon n \g{\n H}{T}+\Theta (S-\overline{S}+\varepsilon (\overline{\Ric}(N,N)-n^2H^2))-n(\varepsilon H\phi+\frac{\partial \phi}{\partial N}).
\end{eqnarray*}
\end{proof}

We can give now the following integral formula for compact hypersurfaces in \N.
\begin{theorem}
\label{thintegral} Let $\Sigma^n$ be a compact hypersurface immersed into \N\ with the general assumptions stated at the beginning of this section. Then \beq \label{integralaux} \int_{\Sigma}
\Theta(S-\overline{S}+\varepsilon\overline{\Ric}(N,N))d\s=n\int_{\Sigma}\frac{\partial \phi}{\partial N}d\Sigma-n(n-1)\varepsilon\int_{\Sigma}H\phi d\Sigma. \eeq
\end{theorem}
\begin{proof}
From \rf{eq2}, and from the fact that $T \in \mathfrak{X}(\N)$ is a conformal Killing field it is immediate to see that
\[
\div(T^\top)=n\phi+nH\Theta.
\]
Therefore, \beq \label{div} \div(HT^\top)=H\div(T^\top)+\g{\n H}{T^\top}=nH\phi+nH^2\Theta+\g{\n H}{T^\top}. \eeq Now, from \rf{laplatheta} and \rf{div} we get
\[
\Delta \Theta+\varepsilon n \div(HT^\top)=\Theta(S-\overline{S}+\varepsilon\overline{\Ric}(N,N))+n(n-1)\varepsilon H\phi-n\frac{\partial \phi}{\partial N}.
\]
Finally, integrating the last expression over the compact hypersurface $\s^n$ and applying the divergence theorem the result follows.
\end{proof}

Theorem\rl{thintegral} becomes especially simple when the field $T$ is a Killing vector field. In that case, \rf{integralaux} is simplified to \beq \label{integral} \int_{\Sigma}
\Theta(S-\overline{S}+\varepsilon\overline{\Ric}(N,N))d\s=0. \eeq On the other hand, when $T$ is a homothetic Killing vector field with $\phi\neq 0$, we may assume without loss of generality that $\phi=1$ and \rf{integralaux} becomes
\[
\int_{\Sigma} \Theta(S-\overline{S}+\varepsilon\overline{\Ric}(N,N))d\s=-n(n-1)\varepsilon \int_{\Sigma}H d\Sigma.
\]

It is especially interesting the case where $T$ is a Killing field, because it will allow us to give some nice consequences of \rf{integral}. Therefore, from now on $T$ will be assumed to be Killing. As a first application, we can
state the following consequence in the case where \N\ is Einstein.
\begin{proposition}
\label{thkilling} Let \N\ be an Einstein space with non-zero scalar curvature $\overline{S}\neq 0$, and endowed with a Killing field (which is assumed to be timelike when \N\ is Lorentzian).
\begin{itemize}
\item[(i)] If \N\ is Lorentzian and $\overline{S}<0$ (respectively, $\overline{S}>0$) there does not exist any compact spacelike hypersurface satisfying $S\leq \overline{S}$ (respectively, $S\geq \overline{S}$).
\item[(ii)] If \N\ is Riemannian and $\overline{S}<0$ (respectively, $\overline{S}>0$) there does not exist any compact hypersurface satisfying $S\leq \overline{S}$ (respectively, $S\geq \overline{S}$) and having $\Theta\neq 0$.
\end{itemize}
\end{proposition}
\begin{proof}
Let us suppose that $\overline{S}<0$ and that there exists a compact hypersurface $\s^n$ such that $S \leq \overline{S}$ in any of the cases (i) or (ii). \N\ being Einstein yields that
$\overline{\Ric}(N,N)=\varepsilon\frac{\overline{S}}{n+1}$, so \rf{integral} becomes \beq \label{intaux} \int_{\Sigma} \Theta\left(S-\overline{S}+\frac{\overline{S}}{n+1}\right)d\s=0. \eeq Recall that when \N\ is Lorentzian the
function $\Theta$ satisfies $\Theta\leq -1<0$. Then, since in the Riemannian case $\Theta$ is assumed to be non-vanishing, we can assume without loss of generalization that $\Theta$ is a negative function over \s. Therefore,
\[
\Theta\left(S-\overline{S}+\frac{\overline{S}}{n+1}\right) \geq 0 \quad \mathrm{on} \quad \s,
\]
and by \rf{intaux} it must vanish. In particular, it must be $S=\overline{S}$ and $\overline{S}=0$, which contradicts the assumption of the theorem. The proof for the case when $\overline{S}>0$ is analogous.
\end{proof}

\section{Hypersurfaces in a product space}
A particular family of spaces with a Killing vector field is that of product spaces. Let $\mene$ be a connected Riemannian surface and consider the product manifold $\mene \times \R$ endowed with the metric
\[
\g{}{}=\g{}{}_\mathbb{M}+\varepsilon dt^2
\]
being $\varepsilon=\pm 1$. When necessary, we will denote the Lorentzian product by $\mene \times \R_1$ in order to distinguish it to the Riemannian one. Observe that $\partial_t$ is a globally defined Killing field over any product
space $\mene \times \R$, with $\g{\partial_t}{\partial_t}=\varepsilon$. Moreover, $\partial_t$ is a parallel vector field, and $\nablabar\pi_{\mathbb{R}}=\varepsilon\partial_t$, where $\pi_\mathbb{R}$ denotes the projection of $\mene
\times \R$ onto the factor \R. The height function $h$ of an immersed hypersurface $\psi:\s^n \fle \mene \times\R$ is the smooth function $h \in \fs$ defined as the projection of the immersion over the factor \R, that is,
$h=\pi_\mathbb{R} \circ\psi$. Therefore, the gradient of $h$ on \s\ is given by
\[
\n h=(\nablabar\pi_\mathbb{R})^\top=\varepsilon \partial_t^\top,
\]
and from the decomposition $\partial_t=\partial_t^\top+\varepsilon \Theta N$ we immediately get \beq \label{normgradh} \|\n h\|^2=\varepsilon(1-\Theta^2). \eeq On the other hand, since $\partial_t$ is parallel in $\mene \times \R$, we
obtain from \rf{eq2} that \beq \label{hessh} \n_X \n h=\Theta AX \eeq and therefore $\Delta h=\varepsilon n H \Theta$.

Here and in what follows \mene\ is assumed to be Einstein whenever $n\geq 3$. For simplicity, $\kappa$ will denote either the (non-necessarily constant) Gaussian curvature of \mdos\ along \s, when $n=2$, or the constant
$S_\mathbb{M}/n$ when $n\geq 3$, being $S_\mathbb{M}$ the scalar curvature of \mene. Therefore,
\[
\overline{\Ric}(U,U)=\kappa |U^\ast|^2
\]
for any vector field $U$ in $\mathfrak{X}(\mene \times \R)$, where $U^\ast$ stands for the projection of the vector field $U$ onto the factor $\mathbb{M}$, that is, $U=U^\ast+\varepsilon \g{U}{\partial_t}\partial_t$, and $|U^\ast|$
denotes its norm with respect to the original metric $\g{}{}_\mathbb{M}$. Then,
\[
\overline{\Ric}(N,N)=\kappa |N^\ast|^2=\kappa \varepsilon (1-\Theta^2),
\]
and
\[
\overline{S}=n \kappa.
\]
Consequently, when \N\ is a product space $\mene \times \R$, \rf{eqgaussaux} becomes \beq \label{eqgaussprod} S=(n-2)\kappa+2\kappa\Theta^2+2\varepsilon\sum_{i<j}\kappa_i\kappa_j. \eeq On the other hand, the integral \rf{integral} is
written as \beq \label{integralproduct} \int_\Sigma\Theta((S-n\kappa)+\kappa(1-\Theta^2))=0. \eeq

\subsection{Hypersurfaces in a Riemannian product space}\label{SS41}
Let us first consider the case where $\s^n$ is a hypersurface immersed into a Riemannian product space $\mene\times\R$. Recall that we are assuming that \s\ is two-sided and the function $\Theta=\g{N}{\partial_t}$ is globally defined.
We will refer to $\Theta$ as the angle function. If \s\ is locally a graph over \mene\ (that it, transversal to $\partial_t$) then  either $\Theta<0$ or $\Theta>0$ over \s. Thus, the assumption \textit{the angle function does not
change sign is on} \s\ is weaker than \s\ \textit{is a local graph}. As already observed by other authors (see for instance \cite{AD,ADR} and \cite{ER}), the angle function is a natural function to consider if one wants to conclude that a hypersurface is necessarily
a slice since, in that case, by \rf{normgradh} we must have that $\Theta^2=1$.

As a first application of our integral formula \rf{integralproduct}, we get the following result for surfaces. Observe that we are not assuming that neither the Gaussian curvature of \mdos\ nor the Gaussian curvature of $\s^2$ are
constant.
\begin{theorem}
\label{th4} Let \mdos\ be a compact Riemannian surface with non-negative Gaussian curvature, $K_{\mathbb{M}}\geq 0$ (respectively, non-positive Gaussian curvature $K_{\mathbb{M}}\leq 0$), and assume that $K_{\mathbb{M}}>0$
(respectively, $K_{\mathbb{M}}<0$) on a dense subset of \mdos. The only compact two-sided surfaces \s\ in the Riemannian product $\mdos\times\mathbb{R}$ such that its angle function $\Theta$ does not change sign and satisfying $K\geq
K_{\mathbb{M}}$ along $\Sigma$ (respectively, $K\leq K_{\mathbb{M}}$) are the slices.
\end{theorem}
\begin{proof}
We can assume without loss of generality that the angle function satisfies $\Theta\leq 0$. Let us consider the case where $K_{\mathbb{M}}\geq 0$ (the case $K_{\mathbb{M}}\leq 0$ is similar). From our hypothesis we have that
$K\geq\kappa\geq 0$ on $\s$. Hence, since $1-\Theta^2\geq 0$, we get that
\[
\Theta(2(K-\kappa)+\kappa(1-\Theta^2))\leq 0,
\]
and it vanishes at a point $p\in\Sigma$ if and only if
\begin{itemize}
\item[(i)] either $\Theta(p)=0$,
\item[(ii)] or $\Theta(p)<0$, $K(p)=\kappa(p)$ and $\kappa(p)(1-\Theta^2(p))=0$.
\end{itemize}
On the other hand, by \rf{integralproduct} we have
\[
\int_\Sigma\Theta(2(K-\kappa)+\kappa(1-\Theta^2))=0,
\]
and then $\Theta(2(K-\kappa)+\kappa(1-\Theta^2))=0$ at every point. This implies that at every $p\in\Sigma$, either $\Theta(p)=0$ or $\Theta(p)=-1$ (see below). Therefore, by continuity the function $\Theta$ must be constant on
$\Sigma$, with either $\Theta=0$ or $\Theta=-1$. But observe that the case $\Theta=0$ cannot happen, since, from \rf{normgradh} we know that there always exist points at which $\Theta=-1$, the critical points of $h$. It follows that
$\Theta=-1$ on $\Sigma$ which means that the surface is a slice.

It remains to prove that $\Theta(2(K-\kappa)+\kappa(1-\Theta^2))=0$ on $\Sigma$ implies that either $\Theta(p)=0$ or $\Theta(p)=-1$  at every $p\in\Sigma$. When $\Theta(p)=0$ there is nothing to prove. Then, let us assume that
$\Theta(p)<0$. If $\kappa(p)\neq 0$, by (ii) above we have that $\Theta(p)=-1$. On the other hand, if $\kappa(p)=0$ we reason as follows. Since $\Theta(p)<0$, we know that $\Sigma$ is locally a graph around  $p$, that is, there exists
an open subset $\Omega\subset\mathbb{M}^2$ and a smooth function $u:\Omega\fle\R$ such that $W=\{ (x,u(x)) : x\in\Omega \}\subset\Sigma$ is an open neighborhood of $p$, with $p=(x_0,u(x_0))$ for a certain $x_0\in\Omega$. In
particular, $\Theta<0$ on $W$. Recall now that we are assuming that $K_{\mathbb{M}}>0$ on a dense subset of \mdos. Therefore, we may find a sequence $\{x_k\}_{k=1}^{\infty}\subset\Omega$ converging to $x_0$ such that
$K_\mathbb{M}(x_k)>0$ at every $k\geq 1$. Let $p_k=(x_k,u(x_k))\in W$, and observe that $\kappa(p_k)=K_\mathbb{M}(x_k)>0$ and $\Theta(p_k)<0$ for every $k\geq 1$. Then, by (ii) above we get that $\Theta(p_k)=-1$ for every $k\geq 1$,
and by continuity $\Theta(p)=-1$. This finishes the proof.
\end{proof}

For the general $n$-dimensional case, we have the following.
\begin{theorem}
\label{th5} Let \mene\ be a compact Einstein Riemannian manifold, $n\geq 3$, with positive scalar curvature $S_{\mathbb{M}}>0$ (respectively, negative scalar curvature $S_{\mathbb{M}}<0$). The only compact two-sided hypersurfaces \s\
in the Riemannian product $\mene\times\mathbb{R}$ such that its angle function $\Theta$ does not change sign and satisfying $\inf_\Sigma S\geq S_{\mathbb{M}}$ (respectively, $\sup_\Sigma S\leq S_{\mathbb{M}}$) are the slices.
\end{theorem}
\begin{proof}
We use the same ideas of the previous proof, taking into account that the scalar curvature $S_\mathbb{M}$ is necessarily constant. We may assume again that $\Theta\leq 0$ and consider the case where $S_{\mathbb{M}}> 0$ is a
positive constant. From our hypothesis we have that $S\geq n\kappa>0$ on \s, so that
\[
\Theta((S-n\kappa)+\kappa(1-\Theta^2))\leq 0,
\]
and by \rf{integralproduct} it vanishes at every point $p\in\Sigma$. In this case, since $\kappa$ is a positive constant we easily have that $\Theta(p)((S(p)-n\kappa)+\kappa(1-\Theta^2(p)))=0$ at a point $p$ if and only if either
$\Theta(p)=0$ or $\Theta(p)=-1$. The same argument as in the case $n=2$ implies then that $\Theta=-1$ and $\Sigma$ is a slice.
\end{proof}

As is well known, the only complete Riemannian surfaces with positive constant Gaussian curvature are, up to a homothety, the sphere $\mathbb{S}^2$ and the real projective plane $\mathbb{RP}^2$. Therefore, we have the following
consequence of Theorem\rl{th4}
\begin{corollary}
\label{corriem} Let \mdos\ be either the sphere $\mathbb{S}^2$ or the real projective plane $\mathbb{RP}^2$.
\begin{itemize}
\item[(i)]  The only compact two-sided surfaces \s\ in the Riemannian product $\mdos\times\mathbb{R}$ such that its angle function $\Theta$ does not change sign and having constant Gaussian curvature are the slices.
\item[(ii)]  The only complete two-sided surfaces \s\ in the Riemannian product $\mdos\times\mathbb{R}$ such that its angle function $\Theta$ is bounded away from 0 and having constant Gaussian curvature are the slices.
\end{itemize}
\end{corollary}

Observe that a cylinder $\gamma\times\mathbb{R}$ over a complete curve $\gamma$ in \mdos\ is an example of a complete and non-compact two-sided surface having $\Theta=0$ and constant Gaussian curvature $K=0$. This
shows that our hypothesis on $\Theta$ in (ii) is necessary. In \cite{AEG1} Aledo, Espinar and G\'{a}lvez posed the following question: are the slices the only complete surfaces in $\mathbb{S}^2 \times \R$ with constant Gaussian
curvature $K$ such that $0<K\leq 1$? Later, in \cite{AEG2}, the same authors proved that there do not exist such surfaces when $0<K<1$, but the question remains open for the case $K=1$. Corollary\rl{corriem} answers partially that
question in the following sense: the only complete surfaces in $\mathbb{S}^2 \times \R$ with constant Gaussian curvature and such that its angle function does not change sign are the slices. In particular, this is also true for entire
complete graphs over $\mathbb{S}^2 \times \R$.

\begin{proof}[Proof of Corollary\rl{corriem}]
We can always assume that $\Theta\leq 0$. Our first objective is to see that in case (ii) the surface \s\ is necessarily compact, and then case (ii) follows directly from case (i). Our hypothesis on $\Theta$ in (ii) means that
$\Theta(p)\leq -\delta<0$ at every point of \s\ for certain $\delta>0$. We claim that this implies that \s\ is necessarily compact. Actually, consider the projection $\Pi:\s \fle \mdos$ of \s\ onto the factor \mdos. Given $p\in
\Sigma$ and a tangent vector $v\in T_p\s$, it is not difficult to see that \beq \label{des} \g{d\Pi_p(v)}{d\Pi_p(v)}_\mathbb{M} \geq c \g{v}{v} \eeq where $c=\frac{\delta^2}{1+\delta^2}>0$. Since $\g{}{}$ is a complete Riemannian
metric on \s\, the same holds for the homothetic metric $\widetilde{\g{}{}}=c\g{}{}$. Then, by \rf{des} the map
\[
\Pi:(\Sigma,\widetilde{\g{}{}})\fle(\mdos,\g{}{}_\mathbb{M})
\]
is a local diffeomorphism which increases the distance. Hence, by \cite[Chapter VIII, Lemma 8.1]{KN} $\Pi$ is a covering map. Therefore, since the universal covering of \mdos\ is $\mathbb{S}^2$, we conclude that \s\ must also be
compact.

Let us see now how to prove the case (i). By compactness of \s, there exists a point $p_0\in\s$ such that $h(p_0)=\min_\Sigma h$. Then, by \rf{normgradh} $\Theta(p_0)=-1$ and $\n^2 h_{p_0}(v,v) \geq 0$ for any $v\in T_{p_0}\s$, where
$\n^2 h$ stands for the Hessian operator on \s. In particular, if we consider the orthonormal frame $\{e_1,e_2\}$ of the principal directions, from \rf{hessh} we get that
\[
\n^2h_{p_0} (e_i,e_i)=-\kappa_i(p_0)\geq 0.
\]
Therefore, $\det(A)(p_0)\geq 0$ and from the Gauss equation \rf{eqgaussprod} we conclude that
\[
K=K(p_0)=K_\mathbb{M}\Theta^2(p_0)+\det(A)(p_0)\geq K_\mathbb{M}.
\]
The result follows now as a direct consequence of Theorem\rl{th4}.
\end{proof}

For general $n$ we get.
\begin{corollary}
Let \mene\ be a compact Einstein Riemannian manifold, $n\geq 3$, with positive scalar curvature. The only compact two-sided hypersurfaces \s\ in the Riemannian product $\mene\times\mathbb{R}$  whose angle function $\Theta$ does not
change sign and having constant scalar curvature are the slices.
\end{corollary}
\begin{proof}
The proof is similar to the proof of Corollary\rl{corriem}. In fact, taking again $p_0\in\s$ a point where $p_0=\min_\Sigma h$ we have that $\Theta(p_0)=-1$ and $\kappa_i(p_0)\leq 0$ for every $i=1,\ldots,n$. Therefore, from the Gauss
equation \rf{eqgaussprod} we obtain that
\[
S=S(p_0)\geq
n\kappa=S_\mathbb{M}.
\]
The result is now a direct consequence of Theorem\rl{th5}.
\end{proof}

\subsection{Hypersurfaces in a Lorentzian product space} \label{SS42}
Let us consider now the case where $\s^n$ is a spacelike hypersurface immersed into a Lorentzian product space $\mene\times\R_1$. Recall that in this case, $N$ denotes the future-directed Gauss map of \s\ and the function $\Theta$
satisfies $\Theta=\g{N}{\partial_t}\leq -1<0$. In particular, we always have $\Theta<0$, and $\Theta(p)=-1$ if and only if $p$ is a critical point of $h$. Reasoning now as in the Riemannian case, we can state the following results.
\pagebreak
\begin{theorem}
\label{th8} \
\begin{itemize}
\item[(i)] Let \mdos\ be a compact Riemannian surface with non-negative Gaussian curvature, $K_{\mathbb{M}}\geq 0$ (respectively, non-positive Gaussian curvature $K_{\mathbb{M}}\leq 0$), and assume that $K_{\mathbb{M}}>0$
(respectively, $K_{\mathbb{M}}<0$) on a dense subset of \mdos. The only compact spacelike surfaces \s\ in the Lorentzian product $\mdos\times\mathbb{R}_1$ satisfying $K\leq K_{\mathbb{M}}$ (respectively, $K\geq K_{\mathbb{M}}$) are
the slices.
\item[(ii)] Let \mene\ be a compact Einstein Riemannian manifold, $n\geq 3$, with positive scalar curvature $S_{\mathbb{M}}>0$ (respectively, negative scalar curvature $S_{\mathbb{M}}<0$). The only compact spacelike hypersurfaces \s\ in the
Lorentzian product $\mene\times\mathbb{R}_1$ satisfying $\sup_\Sigma S\leq S_{\mathbb{M}}$ (respectively, $\inf_\Sigma S\geq S_{\mathbb{M}}$) are the slices.
\end{itemize}
\end{theorem}

\begin{corollary}
\label{corlor} \
\begin{itemize}
\item[(i)] Let \mdos\ be either the sphere $\mathbb{S}^2$ or the real projective plane $\mathbb{RP}^2$. The only complete spacelike surfaces \s\ in the Lorentzian product $\mdos\times\mathbb{R}_1$ with constant
Gaussian curvature are the slices.
\item[(ii)] Let \mene\ be a compact Einstein Riemannian manifold, $n\geq 3$, with positive scalar curvature. The only compact spacelike hypersurfaces \s\ in the
Lorentzian product $\mene\times\mathbb{R}_1$  with constant scalar curvature are the slices.
\end{itemize}
\end{corollary}

\section{Entire graphs in a product space}
\label{S5} As a nice consequence of the results stated in the previous section, in the $n=2$ dimensional case is the possibility of stating them in terms of entire graphs in a product space $\mdos \times \R$. Let $\Omega\subseteq\mdos$ be a
connected domain. Every smooth function $u\in\mathcal{C}^\infty(\Omega)$ determines a graph over $\Omega$ given by $\Sigma(u)=\{ (x,u(x)) : x\in\Omega \}\subset\mdos\times\R$. The metric induced on $\Omega$ from the metric on the
ambient space via $\Sigma(u)$ is given by \beq \label{gu} \g{}{}=\g{}{}_\mathbb{M}+\varepsilon du^2. \eeq In particular, in the case when \mdos\ is Lorentzian $\Sigma(u)$ is a spacelike surface in $\mdos \times \R_1$ if and only if
$|Du|^2<1$ everywhere on $\Omega$, where $Du$ denotes the gradient of $u$ with respect to the metric $\g{}{}_\mathbb{M}$ on $\Omega$ and $|Du|^2=\g{Du}{Du}_\mathbb{M}$. A graph is said to be entire if $\Omega=\mdos$.

Let $\s(u)$ be a graph over a domain $\Omega$, which is assumed to be spacelike when the ambient space is Lorentzian, and let us orient it by the normal field given by
\[
N=\frac{1}{\sqrt{1+\varepsilon|Du|^2}}(-\varepsilon\partial_t+Du),
\]
so that
\[
\Theta=\frac{-1}{\sqrt{1+\varepsilon|Du|^2}}<0.
\]
From the Gauss equation \rf{eqgaussprod}, we know that
\[
K=K_\mathbb{M}\Theta^2+\varepsilon\det(A)=\frac{K_\mathbb{M}}{1+\varepsilon|Du|^2}+\varepsilon\det(A).
\]
Take a local orthonormal frame $\{E_1, E_2\}$ on $\Omega$ (with respect to the metric $\g{}{}_\mathbb{M}$) and observe that
\[
\det(A)=\frac{\det(h_{ij})}{\det(g_{ij})}
\]
where, by \rf{gu}, $g_{ij}=\g{E_i}{E_j}=\delta_{ij}+\varepsilon E_i(u)E_j(u)$, and, by a straightforward computation,
\[
h_{ij}=\g{AE_i}{E_j}=-\g{\nablabar_{E_i}N}{E_j}=\frac{-D^2u(E_i,E_j)}{\sqrt{1+\varepsilon|Du|^2}}.
\]
Therefore, the Gaussian curvature of $\Sigma(u)$ is given by
\[
K=\frac{K_\mathbb{M}}{1+\varepsilon|Du|^2}+\varepsilon\frac{\det(D^2u)}{(1+\varepsilon|Du|^2)^2},
\]
with the restriction $|Du|^2<1$, when $\varepsilon=-1$.

Taking this into account, Corollaries\rl{corriem} and\rl{corlor} can be formulated in terms of graphs in the following way.
\begin{corollary}
\label{coro10} Let \mdos\ be either the sphere $\mathbb{S}^2$ or the real projective plane $\mathbb{RP}^2$, and let $K$ be a real constant. The only entire solutions on \mdos\ to either
\[
(1+|Du|^2)^2K=1+|Du|^2+\det(D^2u)
\]
or
\[
(1-|Du|^2)^2K=1-|Du|^2-\det(D^2u), \quad |Du|^2<1,
\]
where $Du$ and $D^2u$ stand for the gradient and the Hessian of a function $u$ on \mdos, are the constant functions, and $K=1$ necessarily. In particular, if $K\neq 1$, there exist no entire solutions on \mdos\ to the two equations
above.
\end{corollary}
More generally, our Theorems\rl{th4} and\rl{th8} allow us to state the following result, where we do not assume that neither the Gaussian curvature of \mdos\ nor the Gaussian curvature of the graph are constant.
\begin{corollary}
\label{coro11} Let \mdos\ be a compact Riemannian surface with non-negative Gaussian curvature, $K_{\mathbb{M}}\geq 0$ (respectively, non-positive Gaussian curvature $K_{\mathbb{M}}\leq 0$), and assume that $K_{\mathbb{M}}>0$
(respectively, $K_{\mathbb{M}}<0$) on a dense subset of \mdos.
\begin{itemize}
\item[(i)] The only entire solutions on \mdos\ to
\[
(1+|Du|^2)^2K=(1+|Du|^2)K_\mathbb{M}+\det(D^2u),
\]
where $K$ is a smooth function on \mdos\ satisfying $K\geq K_\mathbb{M}$ (respectively, $K\leq K_{\mathbb{M}}$) are the constant functions, and $K=K_\mathbb{M}$ necessarily.
\item[(ii)] The only entire solutions on \mdos\ to
\[
(1-|Du|^2)^2K=(1-|Du|^2)K_\mathbb{M}-\det(D^2u), \quad |Du|^2<1,
\]
where $K$ is a smooth function on \mdos\ satisfying $K\leq K_\mathbb{M}$ (respectively, $K\geq K_{\mathbb{M}}$) are the constant functions, and $K=K_\mathbb{M}$ necessarily.
\end{itemize}
\end{corollary}

In contrast to Corollary\rl{coro10}, when $\mdos=\mathbb{H}^2$ is the hyperbolic plane, there exist entire non-trivial solutions to the corresponding constant Gaussian curvature equation, \beq \label{grafo1}
K=\frac{-1}{1+\varepsilon|Du|^2}+\varepsilon\frac{\det(D^2u)}{(1+\varepsilon|Du|^2)^2}, \eeq with the restriction $|Du|^2<1$, when $\varepsilon=-1$. To see that, it will be appropriate to use the Minkowskian model of the hyperbolic
plane. Write $\mathbb{R}^{3}_1$ for $\mathbb{R}^{3}$, with canonical coordinates $x=(x_0,x_1,x_2)$, endowed with the Lorentzian metric \beq \label{lor} \g{}{}=-dx_0^2+dx_1^2+dx_2^2. \eeq The hyperbolic plane $\mathbb{H}^2$ is the
complete simply connected Riemannian surface with sectional curvature $-1$, which is realized as the hyperboloid
\[
\mathbb{H}^2=\{ x\in{\mathbb R}^{3}_1 : \g{x}{x}=-1, x_0\geq 1 \}\subset\mathbb{R}^3_1
\]
endowed with the Riemannian metric induced from ${\mathbb R}^{3}_1$.

\begin{example}
\label{exam}
\begin{rm}
Let us look for non-trivial solutions on $\mathbb{H}^2$ of \rf{grafo1} of the type $u(x)=f(x_0)$ for a certain smooth function $f(x_0)$ with $x_0\geq 1$. For such a function $u$, its gradient is $Du(x)=-f'(x_0)e_0^\top$,
where $e_0^\top$ is the tangent part of $e_0=(1,0,0)$ along $\mathbb{H}^2$, that is, \beq \label{grafo2} e_0=e_0^\top+x_0x. \eeq In particular, $|Du(x)|^2=f'(x_0)^2(x_0^2-1)$. On the other hand, taking covariant derivatives in
\rf{grafo2} we get
\[
D_X(Du)=f''(x_0)\g{X}{e_0^\top}e_0^\top+f'(x_0)x_0X
\]
for every vector field $X$ tangent to $\mathbb{H}^2$. It follows that
\[
\det(D^2u)=x_0f'(x_0)f''(x_0)(x_0^2-1)+x_0^2f'(x_0)^2.
\]
Hence, equation \rf{grafo1} becomes \beq \label{grafo3} (1+\varepsilon f'(x_0)^2(x_0^2-1))^2K=-1-\varepsilon f'(x_0)^2(x_0^2-1)+\varepsilon(x_0f'(x_0)f''(x_0)(x_0^2-1)+x_0^2f'(x_0)^2). \eeq It can be easily checked that the solution
to the differential equation \rf{grafo3} is given, up to a constant, by
\[
f(x_0)=\sqrt{\frac{\varepsilon(1+K)}{-K}}\, \mathrm{log}\left(\sqrt{1-K(x_0^2-1)}+\sqrt{-K}x_0\right)
\]
where $K<0$ and $\varepsilon(1+K)>0$. Specifically,
\begin{itemize}
\item[$\bullet$] when $\varepsilon=1$, $-1<K<0$, and
\item[$\bullet$] when $\varepsilon=-1$, $K<-1$.
\end{itemize}
\end{rm}
\end{example}

Let us observe that when $\varepsilon=1$, any entire graph in $\mathbb{H}^2\times\mathbb{R}$ is complete. In particular, the entire graphs in our Example\rl{exam} are complete. However, when $\varepsilon=-1$ an entire spacelike
graph in $\mathbb{H}^2\times\mathbb{R}_1$ is not necessarily complete (see, for instance, \cite[Example 3.3]{Al}). Even more, there exist entire spacelike graphs with constant Gaussian curvature which are not complete; for instance
\cite[Example 3.1]{Al} is an example of a maximal entire graph with constant Gaussian curvature $-1$ in $\mathbb{H}^2\times\mathbb{R}$ which is not complete.

In our case, the completeness of the entire spacelike graphs given in Example\rl{exam} follows from the more general technical result.
\begin{lemma}
Let $\Sigma(u)$ be an entire spacelike graph in $\mathbb{M}^2\times\mathbb{R}_1$, where \mdos\ is a complete Riemannian surface. If
\[
\sup_{x\in\mathbb{M}}|Du(x)|^2<1,
\]
where $|Du|^2=\g{Du}{Du}_\mathbb{M}$, then $\Sigma(u)$ is complete.
\end{lemma}
\begin{proof}
Observe that the induced metric on $\mathbb{M}^2$ via the graph $\Sigma(u)$ is given by
\[
g(v,v)=\g{v}{v}_{\mathbb{M}}-\g{v}{Du}^2_{\mathbb{M}},
\]
and by Cauchy-Schwarz inequality we have
\[
g(v,v)\geq \g{v}{v}_{\mathbb{M}}(1-|Du|^2)\geq \g{v}{v}_{\mathbb{M}}(1-\sup_{x\in\mathbb{M}}|Du(x)|^2).
\]
Finally, since $(1-\sup_{x\in\mathbb{M}}|Du(x)|^2)>0$ and $\g{}{}_\mathbb{M}$ is complete, the metric $g$ is also complete on $\mathbb{M}^2$.
\end{proof}
Observe that in our case, for every $x\in\mathbb{H}^2$
\[
|Du(x)|^2=f'(x_0)^2(x^2_0-1)=\frac{-(1+K)(x^2_0-1)}{1-K(x^2_0-1)}, \quad \mathrm{with} \quad K<-1,
\]
and then
\[
\sup_{x\in\mathbb{H}^2}|Du(x)|^2=\sup_{x_0\geq 1}f'(x_0)^2(x_0^2-1)=1+\frac{1}{K}<1.
\]

\section*{Acknowledgements}
The authors would like to thank to the referee for valuable suggestions which improved the paper.

\end{document}